
\input amstex
\NoBlackBoxes
\documentstyle{amsppt}
\pagewidth{125mm}
\pageheight{185mm}
\pageno=1

\topmatter
\title Resonances and Eigenvalues for the Constant Mean Curvature Equation
\endtitle
\dedicatory{To Carlos Kenig in friendship and admiration}\enddedicatory
\author{Sagun Chanillo}\endauthor

\address RUTGERS UNIVERSITY, DEPT. OF MATH., HILL CENTER, BUSCH CAMPUS\endgraf
 110 FRELINGHUYSEN RD, PISCATAWAY, NJ 08854, USA
\endaddress
\email
chanillo\@math.rutgers.edu
\endemail
\rightheadtext{Resonances and Eigenvalues}
\topmatter
\abstract
In this paper we study resonances and eigenvalues for the nonlinear constant mean curvature eqn. linearized around the bubbles found by Brezis-Coron. This nonlinear eqn. is also called a H-system eqn. For degree one bubbles we only find resonances. For higher degree we prove eigenvalues occur. Our goal is to eventually obtain dispersive estimates for the wave eqn. associated to the linear and non-linear problem, a study of which was initiated by Chanillo-Yung.

\endabstract
\endtopmatter
\document
\bigskip
\subhead 1. Introduction\endsubhead
\medskip
Let
$$
u: \Omega \subset \Bbb R^2 \to \Bbb R^3,
$$
where $\Omega$ is a domain.  We will denote vector cross products of two vectors $u, v$ by $u \wedge v$.  Consider now the equ.

$$
\Delta u = 2 u_x \wedge u_y , u= (u_1, u_2, u_3)
\tag 1.1
$$
where $\Delta u= (\Delta u_1, \Delta u_2, \Delta u_3)$.
\medskip

If $u$ satisfies (1.1) and in addition the Plateau conditions
$$
|u_x|^2 = |u_y|^2, u_x \cdot u_y = 0,
\tag 1.2
$$
then the range of $u$ defines a constant mean curvature surface in $\Bbb R^3$.  In this work we shall ignore the Plateau conditions (1.2).  (1.1) can be written in variational form
$$
\Cal E (u) = \frac{1}{2} \int_\Omega |\nabla u |^2 + \frac{2}{3} \int_\Omega u \cdot u_x \wedge  u_y.
\tag 1.3
$$
Brezis-Coron [1] and Struwe [10] studied (1.3) and its Morse theory.  Brezis-Coron [2] classified all the finite energy solutions of (1.1).  These finite energy solutions will be referred to as bubbles in the sequel.  [2] also obtained the first result for (1.1) as to how compactness in (1.3) fails.  (1.1) is a conformally invariant PDE and also invariant by translation and rotation in the target.  Chanillo-Malchiodi [4] studied the Morse theory further for (1.1) and constructed multi-bubble solutions to (1.1).  We point out that the Topology of $\Omega$ plays a strong role in constructing solutions.  For example, a result of Wente [11], another proof may be found in Prop 3.1 of [4], shows that if $\Omega$ is simply connected, then
$$
\Delta u = 2 u_x \wedge u_y,  u\Bigg|_{\partial \Omega} \equiv 0
\tag 1.4
$$
has only the solution $u\equiv 0$, while if $\Omega$ is an annulus, Wente constructed non-trivial solutions to (1.4) in [11].

We now recall the classification result of [2].  Consider those solutions in all $\Bbb R^2$ satisfying,
$$
\Delta u = 2 u_x \wedge u_y, \int_{\Bbb R^2} |\nabla u |^2 < \infty.
$$
Then, we have quantization, that is necessarily
$$
\int_{\Bbb R^2} |\nabla u |^2 = 8\pi m, m \in \Bbb N, m \geq 1.
$$
and,
$$
u = \pi^{-1} \left( \frac{P(z)}{Q(z)}\right), \pi : S^2\to \Bbb R^2, z = (x, y).
$$
$\pi$ is the stereographic projection to the plane from the Riemann sphere, $P(z), Q(z)$ are polynomials which are holomorphic and $m = \max\{\deg P, \deg Q\}$.  Thus the degree 1 bubbles which are basic have the form
$$
(x, y) \to \left(\frac{2x}{1+|z|^2} , \frac{2y}{1+|z|^2}, \frac{|z|^2-1}{|z|^2 +1}\right)
\tag 1.5
$$
where $z= x+iy$, which we will abuse sometimes and think $z=(x,y)$.  Degree $m$ bubbles, $m\geq 2$, can be written as (note these are special degree $m$ bubbles)
$$
\left( \frac{2 \Re \, z^m}{1+ |z|^{2m}},\,  \frac{ 2 \Im \, z^m}{1+|z|^m},\,  \frac{|z|^{2m} - 1}{|z|^{2m} + 1}\right).
\tag 1.6
$$

In [5] the study of the wave eqn. corresponding to (1.1)
$$
\Delta u - u_{tt} = 2u_x\wedge u_y,
u \Bigg|_{t=0} = u_0 , u_t \Bigg|_{t=0}  = u_1
\tag 1.7
$$
was initiated.  The right side of (1.7) is an example of a null form.  Since we are in dimension 2, with poor Strichartz estimates and lack of Huygens principle, many problems concerning (1.7) remain unresolved.  Some attempts have been made [3] to understand random data versions of (1.7).  One of the main results of [5], in the spirit of a similar result for another conformally invariant eqn., the Yamabe eqn. [8], is that if
$$
\| \nabla u_0 \|^2_2 > 8\pi
$$
which is saying that if the initial data has more energy than a degree 1 bubble, and with an additional condition on $u_1$, involving energy trapping, then the solution to (1.7) blows up in finite time.  Furthermore in [6], via unique continuation arguments it is proved that the blow up is not self-similar.  No refined bubbling analysis for (1.7) is known to date.

Our goal here is to take an initial step for a linearized version of (1.7) and study dispersion and scattering like what was done in 3D for the Yamabe eqn. in [9].  To do so one has to study Born expansions.  To perform such expansions one needs to know if resonances and eigenvalues exist for the linearized elliptic part.  Appearance of eigenvalues in the spectrum complicates matters.

To state our theorems, we first linearize (1.1) around a bubble $u$, and get
$$
\Delta w = 2w_x\wedge u_y + 2u_x \wedge w_y.
\tag 1.8
$$
We now prefer to study (1.8) on the punctured sphere $S^2\setminus\{N\}$, where $N$ is the North Pole.  So we study (1.8) on $\Bbb R^2$, equipped with the metric
$$
\frac {4}{(1+|z|^2)^2} (dx^2 + dy^2), z= (x,y).
$$
(1.8) becomes due to a conformal change of the metric

$$
\frac{1}{\varphi} \Delta w = \frac {2}{\varphi} ( w_x \wedge u_y + u_x \wedge w_y), \varphi = \frac {4}{(1+|z|^2)^2}.
\tag 1.9
$$
Thus the eigenvalue eqn. for (1.9) is
$$
\Delta w = 2 ( w_x \wedge u_y + u_x \wedge w_y) + \frac{4\lambda}{(1+|z|^2)^2} w
\tag 1.10
$$
where $z=(x, y)$ and
$$
\Delta w = \frac {\partial^2 w}{\partial x^2} + \frac{\partial^2 w}{\partial y^2}.
$$
The eqn. (1.10) when $\lambda = 0$ is of course
$$
\Delta w = 2(w_x \wedge u_y + u_x \wedge w_y).
\tag 1.11
$$
We also notice by integration by parts the self-adjoint property
$$
\Cal E (v,w) = \Cal E(w, v) = \int_{\Bbb R^2} \nabla v \cdot \nabla w + 2\int_{\Bbb R^2} v \cdot (w_x\wedge u_y + u_x \wedge w_y)
\tag 1.12
$$
for $v, w \in C^\infty_0(\Bbb R^2, \Bbb R^3)$.  In view of (1.12) we search for solutions to (1.10), (1.11) for which $w\in L^2 (\Bbb R^2)$.

\proclaim{Definition}

(a) We say $\lambda$ is a resonance for (1.10), (1.11), if $w \not\in L^2 (\Bbb R^2)$ and $w$ satisfies (1.10), (1.11).
\medskip

(b) We say $\lambda$ is an eigenvalue if $w\in L^2(\Bbb R^2)$ and $w$ satisfies (1.10), (1.11).
\endproclaim
Lastly we shall only focus on $w$ that is co-rotational.  That is keeping in mind the expression for the bubbles (1.6), we assume for $(r, \theta)$ polar coordinates,
$$
w(r, \theta) = (f(r)\cos m \theta, f(r) \sin m \theta, g(r)).
\tag 1.13
$$

\proclaim{Theorem 1.1}

(a) $\lambda = 0$ is a resonance for (1.11) for $u$ a degree 1 bubble.

\medskip

(b) $\lambda = 0$ is an eigenvalue for (1.10), (1.11), for $u$ a degree $m$ bubble of the form (1.6), $m\geq 2$, and with an eigenfunction of the form (1.13).
\endproclaim

\proclaim{Theorem 1.2}  For $u$ a degree 1 bubble, (1.10), (1.11) has no eigenfunctions of the form (1.13) for any $\lambda$.  That is there are no co-rotational eigenfunctions for any choice of $\lambda$.
\endproclaim

We end by recording a problem.

\medskip

\noindent{\bf Problem}.  In view of Theorem 1.1, when $u$ is a degree $m$ bubble, $m\geq 2$ we know $\lambda = 0$ is an eigenvalue.  Does (1.10) have other eigenvalues besides $\lambda = 0$, where the eigenfunctions are co-rotational of the form (1.13), when $u$ is a degree $m$ bubble given by (1.6), $m\geq 2$?


\subhead 2. Proofs of the Results \endsubhead


We begin by proving Theorem (1.1).

\demo{Proof}  First we note that
$$
\Delta u = 2u_x \wedge u_y
$$
is invariant under dilations, that is under the map $z\to \delta z, \delta >0, z = x+iy$.  We only consider dilations as we are focused on co-rotational solutions to our linearized eqn.
$$
\Delta w = 2 u_x \wedge w_y + 2 w_x \wedge u_y +\frac{4\lambda}{(1+|z|^2)^2} w
\tag 2.1
$$
where $u(x, y)$ is a bubble of degree $m$.  Since we are examining eigenvalues and resonances at $\lambda = 0$, we set $\lambda = 0$ in (2.1) and also
$$
u(z) = \left( \frac {2 \Re z^m}{1+|z|^{2m}}, \frac {2 \Im z^m}{1+|z|^{2m}}, 1-\frac{2}{1+|z|^{2m}}\right).
\tag 2.2
$$
Dilation of the bubble given by (2.2) yields
$$
u^\delta(x, y) = \left(\frac{2\delta^m \Re z^m}{1+\delta^{2m}|z|^{2m}}, \,
\frac{2 \delta^m \Im z^m}{1+\delta^{2m}|z|^{2m}}, 1 - \frac{2}{1+\delta^{2m}|z|^{2m}}\right).
\tag 2.3
$$
By invariance we have
$$
\Delta u^\delta = 2 u^\delta_x \wedge u_y^\delta.
\tag 2.4
$$
Differentiating (2.3), (2.4) in $\delta$, and setting $\delta =1$, we get,
$$
w= \frac{\partial u^\delta}{\partial \delta}\Bigg|_{\delta= 1} = \left( \frac{2 m \Re \, z^m (1-|z|^{2m})}{(1+|z|^{2m})^2},  \frac{2 m \Im \, z^m(1-|z|^{2m})}{(1+|z|^{2m})^2}, \frac{2 m |z|^{2m}}{(1+|z|^{2m})^2} \right),
\tag 2.5
$$
and
$$
\Delta w = 2 w_x \wedge u_y + 2 u_x \wedge w_y.
\tag 2.6
$$
We see from (2.5), that there exists $c_1, c_2 > 0$ such that
$$
|w (x, y)|\leq c_1(1+|z|)^{-m}
$$
and
$$
c_2 |z|^{-m} \leq |w(x, y)|, \text{ for } |z| \geq 10.
$$
Thus $w\in L^2 (\Bbb R^2)$, for $m\geq 2$, and $w \not\in L^2(\Bbb R^2)$ when $m=1$.  Hence we conclude that $\lambda = 0$ is an eigenvalue for linearized problem around a bubble of degree $\geq 2$, and $\lambda = 0$ is a resonance for the linearized problem when the CMC eqn. is linearized around a bubble of degree one.  When $\lambda=0$ and $u$ a degree one bubble, we will conclusively rule out eigenvalues in the next theorem. This ends the proof.
\enddemo

We now prove Theorem (1.2).

\demo{Proof}  Our first goal is to compute the ODE for the co-rotational soln. for the linearized problem, where we linearize at a degree $m$ bubble.  We will later specialize to $m=1$.  Since we are considering the punctured sphere $S^2\setminus \{N\}$, where $N$ is the North pole, as explained in the introduction, our eigenvalue eqn. is
$$
\Delta w = 2 w_x \wedge u_y +2 u_x \wedge w_y + \frac{ 4\lambda}{(1+|z|^2)^2} w
\tag 2.7
$$
where $u$ is a degree $m$ bubble. $w$ being co-rotational has the form in polar coordinates $(r,\theta), w(r, \theta) = (f(r) \cos m \theta, f (r) \sin m \theta, g(r))$.  By changing variables it is easily seen
$$
w_x \wedge u_y + u_x\wedge w_y = \frac{1}{r} w_r \wedge u_\theta +\frac{1}{r} u_r \wedge w_\theta.
\tag 2.8
$$
Next we perform a Kelvin transformation by setting

$$
r = \frac{1}{t}, \theta \to -\theta
$$
and routine computation yields that (2.8) becomes
$$
t^3(w_t \wedge u_\theta + u_t \wedge w_\theta).
\tag 2.9
$$
The Laplacian transformes to
$$
t^4\left( \frac{\partial^2}{\partial t^2} + \frac{1}{t} \frac{\partial}{\partial t} + \frac{1}{t^2} \frac{\partial^2}{\partial\theta^2}\right) = t^4\Delta.
$$
After performing the Kelvin transformation
$$
\tilde{w} (t, \theta) = (\tilde{f}(t) \cos m \theta, - \tilde{f}(t)\sin m \theta, \tilde{g}(t))
$$
where
$$
\tilde{f}(t) = f(\frac{1}{t}) , \tilde{g}(t) = g(\frac{1}{t}).
$$
Likewise the bubble becomes
$$
\tilde{u} = \left( \tilde{F}(t) \cos m \theta, - \tilde{F}(t) \sin m \theta, \tilde{G}(t)\right).
$$
We conserve notation and drop the tildes.

$$
\aligned
w_t    &=(f'(t))\cos m \theta, - f'(t) \sin m \theta, g'(t))   \\
w_\theta &= (-mf \sin m \theta, - mf \cos m \theta, 0)    \\
u_t     &= (F'\cos m \theta, -F' \sin m \theta, G')\\
u_\theta &=(-m F \sin m \theta, - m F \cos m \theta, 0).
\endaligned
$$
We obtain by elementary computation,
$$
\aligned
w_t \wedge u_\theta &= (m F g'\cos m \theta, - m F g' \sin m \theta, - m Ff')\\
u_t \wedge w_\theta &= ( mfG' \cos m \theta, - m f G' \sin m \theta, - mfF').
\endaligned
$$
Hence
$$
2(w_t \wedge u_\theta + u_t \wedge w_\theta) = 2\left(( fG' + Fg') \cos m \theta, - m(f G' + Fg')\sin m \theta, - m(Ff' + fF')\right).
$$
Hence we arrive at a pair of second order ODE,
$$
t^4 \left(f'' + \frac{1}{t} f' - \frac{m^2}{t^2} f\right) =  2m t^3(fG'+Fg')+
\frac{4\lambda t^4}{(1+t^2)^2} f .
\tag 2.10
$$
$$
t^4 \left(g'' + \frac{1}{t} g'\right) =  -2m t^3(Ff'+ fF')+
\frac{4\lambda t^4}{(1+t^2)^2 } g.
\tag 2.11
$$
Before we proceed further we note that if $u$ is a degree $m$ bubble,
$$
|\nabla u| \leq C.
$$
Thus from (2.9) and (2.7), any eigenfunction will satisfy
$$
\aligned
t^4 |\Delta w| &\leq t^4 |\nabla w| + \frac{4t^4}{(1+t^2)^2} \,  |\lambda | |w|\\
               &\leq C(|\nabla w | + |w|) t^4.
\endaligned
$$
We have
$$
|\Delta w| \leq C(|\nabla w| + |w|).
\tag 2.12
$$
If we show that $w$ vanishes to infinite order at $t=0$, the origin, then by applying the unique continuation Lemma 2.6.1, pg 70 of [7], we may conclude that $w \equiv 0$ in $\Bbb R^2$.  This will complete the proof of the theorem.
To show that $w$ does indeed vanish to infinite order, we need to assume by contradiction that $f, g \in L^2(\Bbb R^2)$ and $f, g$ are of course smooth.  Thus both $f, g$ have a formal power series expansion.

$$
f(t) \sim \sum_{n\geq 2} a_n t^n , \,  g(t) \sim \sum_{n\geq 2} b_n t^n.
$$
Obviously $n \geq 2$, for we are assuming $f, g \in L^2( \Bbb R^2)$; and $t$ is the variable after Kelvin transformaion.  The expression for the $m=1$ bubble in $(t,\theta)$ coordinates is
$$
\left(\frac{2t}{1+t^2} \cos \theta, - \frac{2t}{1+t^2} \sin \theta,  \frac{1-t^2}{1+t^2}
\right).
$$
Thus
$$
F(t) = \frac{2t}{1+t^2}, G(t) =- 1 + \frac{2}{1+t^2}.
$$
Thus,
$$
F'(t) = \frac{2(1-t^2)}{(1+t^2)^2} , G' =- \frac{4t}{(1+t^2)^2}.
$$
Now
$$
Ff' + fF' = \left(\sum_{n\geq 2} n a_n t^{n-1}\right) \frac{2t}{1+t^2} + 2
\left(\sum_{n\geq 2} a_n t^n\right)  \frac{(1-t^2)}{(1+t^2)^2}.
$$
Thus, (2.11) becomes after cancellling off $t^4$ from both sides; for $m=1$,
$$
g'' + \frac{1}{t} g' = - \frac{4}{1+ t^2} \sum_{n\geq 2} n a_n t^{n-1} + 4
\left(\sum_{n\geq 2}  a_n t^{n-1}\right)\frac{(1-t^2)}{(1+t^2)^2} + \frac{4\lambda}{(1+t^2)^2} g.
$$
Since $g,(t) = \sum_{n\geq 2}  b_n t^n$, the left side above,
is
$$
\sum_{n\geq 2} n^2 b_n t^{n-2} = \sum_{n\geq 0} (n+2)^2 b_{n+2} t^n.
$$
We have using (2.11)

$$
\aligned
 &(1+t^2)^2 \sum_{n\geq 0} (n+2)^2 b_{n+2} t^n\\
 &= - 4(1+t^2) \sum_{n\geq 2} n a_n t^{n-1} + 4\left(\sum_{n\geq 2}  a_n t^{n-1}\right)(1-t^2) + 4\lambda \sum_{n\geq 2}  b_n t^n.
\endaligned
\tag 2.13
$$
It is easily seen from (2.13), that
$$
(n+2)^2 b_{n+2} =  \sum_{k\leq n+1} \alpha_k b_k +  \sum_{k\leq n+1} \gamma_k a_k , n\geq 0,
$$
for suitable constants $\alpha_k, \gamma_k$.  We may apply induction to conclude that $b_k = 0$, for all $k\geq 0$.  Thus $g(t)$ vanishes to infinite order at $t=0$.  A similar arguement but now using (2.10), allows us to express
$$
\left((n+2)^2 - 1\right)  a_{n+2} = \sum_{k\leq n+1} \sigma_k b_k + \sum_{k\leq n+1} \delta_k a_k,  n\geq 0.
$$
We conclude $a_k = 0$, for all $k\geq 0$. Thus $f(t)$ vanishes to infinite order at $t=0$.
This then proves our theorem.

\enddemo


\Refs\nofrills{References}

\widestnumber\key{ZAAAAAAA}


\ref
\key 1
\by Brezis H., Coron J.-M.
\paper Multiple solutions of H-systems and Rellich's conjecture
\jour Comm. Pure and Appl. Math.
\vol 37(2)
\yr 1984
\pages 149--187
\endref
\medskip



\ref
\key 2
\bysame
\paper Multiple solutions of H-systems or how to blow bubbles
\jour Arch. Rational Mech. Analysis
\vol 89(1)
\year 1985
\pages 21--56
\endref
\medskip

\ref
\key 3
\by  Chanillo S., Czubak M., Mendelson D., Nahmod A., Staffiliani G.
\paper Almost sure boundedness of iterates for derivative nonlinear wave eqns.
\jour Comm. Analysis and Geom.
\vol 28(4)
\yr 2020
\pages 941--975
\endref
\medskip


\ref
\key 4
\by  Chanillo S., Malchiodi A.
\paper  Asymptotic Morse theory for $\Delta v = 2v_x\wedge v_y$
\jour Comm. Analysis and Geom.
\vol 13(1)
\yr 2005
\pages 187--251
\endref
\medskip

\ref
\key 5
\by  Chanillo S., Yung P.-L.
\paper Wave Equations associated to Liouville systems and constant mean curvature eqns.
\jour Advances in Math.
\vol 235
\yr 2013
\pages 187--207
\endref
\medskip


\ref
\key 6
\bysame
\paper Absence of self-similar blow-up and local well-posedness of the constant mean curvature wave eqn.
\jour J. Funct. Analysis
\vol 269
\yr 2015
\pages 1180--1202
\endref
\medskip

\ref
\key 7
\by Jost J.
\book Two-dimensional geometric variational problems
\publ Wiley Interscience Publication
\vol
\yr 1991
\pages
\endref
\medskip




\ref
\key 8
\by Kenig C., Merle F.
\paper Global well-posedness, scattering and blow-up for the energy critical focusing wave eqn.
\jour Acta Math.
\vol 201(2)
\yr 2008
\pages  147--212
\endref
\medskip

\ref
\key 9
\by Krieger J., Schlag W.
\paper On the focusing critical semi-linear wave eqn.
\jour Amer. J. of Math
\vol 129(3)
\yr 2007
\pages 843--913
\endref
\medskip



\ref
\key 10
\by  Struwe M.
\paper\nofrills Plateau's problem and the calculus of variations
\inbook Mathematical Notes
\publ Princeton Univ. Press, Princeton NJ
\vol 35
\yr 1988
\pages x + 148
\endref
\medskip



\ref
\key 11
\by Wente H.
\paper  The differential eqn., $\Delta x = 2H x_u \wedge x_v$ with vanishing boundary values
\jour Proc Amer. Math. Soc.
\vol 50
\yr 1975
\pages 131--137
\endref
\medskip







\endRefs
\enddocument